\documentclass[a4paper,11pt,reqno]{amsart}

\newtheorem{theorem}[equation]{Theorem}
\newtheorem{lemma}[equation]{Lemma}

\newtheorem{proposition}[equation]{Proposition}

\numberwithin{equation}{section}

\theoremstyle{definition}

\newtheorem*{example*}{Example}

\newtheorem{remark}[equation]{Remark}
\newtheorem*{remark*}{Remark}

\newcommand{\bN}{{\mathbb N}}
\newcommand{\bZ}{{\mathbb Z}}

\newcommand{\cJ}{{\mathcal J}}
\newcommand{\cC}{{\mathcal C}}
\newcommand{\cT}{{\mathcal T}}
\newcommand{\frg}{{\mathfrak g}}

\newcommand{\frh}{{\mathfrak h}}

\newcommand{\subo}{_{\bar 0}}
\newcommand{\subuno}{_{\bar 1}}

\providecommand{\espan}[1]{\text{span}\left\{ #1\right\}}

 \DeclareMathOperator{\frsl}{{\mathfrak{sl}}}
 \DeclareMathOperator{\frsp}{{\mathfrak{sp}}}
 \DeclareMathOperator{\frso}{{\mathfrak{so}}}

 \DeclareMathOperator{\frosp}{{\mathfrak{osp}}}

  \DeclareMathOperator{\tr}{tr}
 \DeclareMathOperator{\ad}{ad}
  \DeclareMathOperator{\Ad}{Ad}
 \DeclareMathOperator{\der}{\mathfrak{der}}
 \DeclareMathOperator{\inder}{\mathfrak{inder}}
 \DeclareMathOperator{\End}{End}
 \DeclareMathOperator{\Hom}{Hom}
 \DeclareMathOperator{\Mat}{Mat}

 \DeclareMathOperator{\charac}{char}
 \DeclareMathOperator{\Cl}{\mathfrak{Cl}}

\newenvironment{romanenumerate}
 {\begin{enumerate}
 \renewcommand{\itemsep}{4pt}
 }{\end{enumerate}}

\begin{document}

\title{Some new simple modular Lie superalgebras}

\author[Alberto Elduque]{Alberto Elduque$^{\star}$}
 \thanks{$^{\star}$ Supported by the Spanish Ministerio de
 Educaci\'{o}n y Ciencia
 and FEDER (MTM 2004-081159-C04-02) and by the
Diputaci\'on General de Arag\'on (Grupo de Investigaci\'on de
\'Algebra)}
 \address{Departamento de Matem\'aticas, Universidad de
Zaragoza, 50009 Zaragoza, Spain}
 \email{elduque@unizar.es}

\date{\today}



\begin{abstract}
Two new simple modular Lie superalgebras will be obtained in
characteristics $3$ and $5$, which share the property that their
even parts are orthogonal Lie algebras and the odd parts their spin
modules. The characteristic $5$ case will be shown to be related, by
means of a construction of Tits, to the exceptional ten dimensional
Jordan superalgebra of Kac
\end{abstract}

\maketitle


\section{Introduction}

There are well-known constructions of the exceptional simple Lie
algebras of type $E_8$ and $F_4$ which go back to Witt
\cite{Witt41}, as $\bZ_2$-graded algebras
$\frg=\frg\subo\oplus\frg\subuno$ with even part the orthogonal Lie
algebras $\frso_{16}$ and $\frso_9$ respectively, and odd part given
by their spin representations (see \cite{Adams}).

Brown \cite{Brown82} found a new simple finite dimensional Lie
algebra over fields of characteristic $3$ which presents the same
pattern, but with $\frg\subo=\frso_7$.

For the simple Lie superalgebras in Kac's classification
\cite{KacLie}, only the orthosymplectic Lie superalgebra
$\frosp(1,4)$ presents the same pattern, since $\frg\subo=\frsp_4$
here, and $\frg\subuno$ is its natural four dimensional module. But
$\frsp_4$ is isomorphic to $\frso_5$ and, as such, $\frg\subuno$ is
its spin module.

Quite recently, the author \cite{Eld05} found another instance of
this phenomenon. There exists a simple Lie superalgebra over fields
of characteristic $3$ with even part isomorphic to $\frso_{12}$ and
odd part its spin module.

This paper is devoted to settle the question of which other simple
either $\bZ_2$-graded Lie algebras or Lie superalgebras present this
same pattern: the even part being an orthogonal Lie algebra and the
odd part its spin module.

It turns out that, besides the previously mentioned examples, and of
$\frso_9$, which is the direct sum of $\frso_8$ and its natural
module, but where, because of triality, this natural module can be
substituted by the spin module, there appear exactly two other
possibilities for Lie superalgebras, one in characteristic $3$ with
even part isomorphic to $\frso_{13}$, and the other in
characteristic $5$, with even part isomorphic to $\frso_{11}$. These
simple Lie superalgebras seem to appear here for the first time.

The characteristic $5$ case will be shown to be strongly related to
the ten dimensional simple exceptional Kac Jordan superalgebra, by
means of a construction due to Tits. As has been proved by McCrimmon
\cite{McC}, and indirectly hinted in \cite{EOpseudo}, the Grassmann
envelope of this Jordan superalgebra satisfies the Cayley-Hamilton
equation of degree $3$ and hence, as shown in \cite{BZ} and
\cite{BETits}, this Jordan superalgebra $\cJ$ can be plugged into
the second component of Tits construction \cite{Tits66}, the first
component being a Cayley algebra. The even part of the resulting Lie
superalgebra is then isomorphic to $\frso_{11}$ and the odd part
turns out to be its spin module.

The characteristic $3$ case is related to the six dimensional
composition superalgebra $B(4,2)$ (see \cite{EOsupercompo} and
\cite{Shestakov}) and, therefore, to the exceptional Jordan
superalgebra of $3\times 3$ hermitian matrices $H_3(B(4,2))$. This
will be left to a forthcoming paper \cite{CE}, where an Extended
Freudenthal Magic Square in characteristic $3$ will be considered.

\smallskip

\emph{Throughout the paper, $k$ will always denote an algebraically
closed field of characteristic $\ne 2$.}

\smallskip

The paper is organized as follows. In the next section, the basic
properties of the orthogonal Lie algebras, associated Clifford
algebras and spin modules will be reviewed in a way suitable for our
purposes.

In Section 3, the simple $\bZ_2$-graded Lie algebras and the simple
Lie superalgebras whose even part is an orthogonal Lie algebra of
type $B$ and its odd part its spin module will be determined. The
two new simple Lie superalgebras mentioned above appear here.
Section 4 is devoted to type $D$, and here the objects that appear
are either classical or a Lie superalgebra in characteristic $3$
with even part $\frso_{12}$, which appeared for the first time in
\cite{Eld05} related to a Freudenthal triple system, which in turn
is constructed in terms of the Jordan algebra of the hermitian
$3\times 3$ matrices over a quaternion algebra.

Finally, Section 5 is devoted to study the relationship of the
exceptional Lie superalgebra that has appeared in characteristic
$5$, with even part isomorphic to $\frso_{11}$, to the Lie
superalgebra obtained by means of Tits construction in terms of the
 Cayley algebra and of the exceptional ten dimensional Jordan
superalgebra of Kac.

\bigskip

\section{Spin modules}

Let $V$ be a vector space of dimension $l\geq 1$ over the field $k$,
let $V^*$ be its dual vector space, and consider the $(2l+1)$
dimensional vector space $W=ku\oplus V\oplus V^*$, with the regular
quadratic form $q$ given by
\begin{equation}\label{eq:q}
q(\alpha u+v+f)=-\alpha^2+f(v),
\end{equation}
for any $\alpha\in k$, $v\in V$ and $f\in V^*$.

Let $\Cl(V\oplus V^*,q)$ be the Clifford algebra of the restriction
of $q$ to $V\oplus V^*$, and let $\Cl\subo(W,q)$ be the even
Clifford algebra of $q$. As a general rule, the multiplication in
Clifford algebras will be denoted by a dot: $x\cdot y$. The linear
map
\[
V\oplus V^*\rightarrow \Cl\subo(W,q):\ x\mapsto u\cdot
x=\frac{1}{2}(u\cdot x-x\cdot u)=\frac{1}{2}[u,x]^{\cdot},
\]
extends to an algebra isomorphism
\begin{equation}\label{eq:Psi}
\Psi: \Cl(V\oplus V^*,q)\rightarrow \Cl\subo(W,q).
\end{equation}
Moreover, let $\tau$ be the involution (that is, involutive
antiautomorphism) of $\Cl(W,q)$ such that $\tau(w)=w$ for any $w\in
W$, and let $\tau\subo$ be its restriction to $\Cl\subo(W,q)$. On
the other hand, let $\tau'$ be the involution of $\Cl(V\oplus
V^*,q)$ such that $\tau'(x)=-x$ for any $x\in V\oplus V^*$. Then,
for any $x\in V\oplus V^*$,
\[
\tau\subo\bigl(\Psi(x)\bigr)=\tau\subo(u\cdot x)=\tau(x)\cdot\tau(u)
 =x\cdot u=-u\cdot x=u\cdot\tau'(x)=\Psi\bigl(\tau'(x)\bigr),
\]
so that $\Psi$ in \eqref{eq:Psi} is actually an isomorphism of
algebras with involution:
\begin{equation}\label{eq:Psitaus}
\Psi:\bigl(\Cl(V\oplus V^*,q),\tau'\bigr)\rightarrow
\bigl(\Cl\subo(W,q),\tau\subo\bigr).
\end{equation}

\medskip

Consider now the exterior algebra $\bigwedge V$. Multiplication here
will be dented by juxtaposition. This conveys a natural grading over
$\bZ_2$: $\bigwedge V=\bigwedge\subo V\oplus\bigwedge\subuno V$. In
other words, like Clifford algebras, $\bigwedge V$ is an associative
superalgebra. For any $f\in V^*$, let $df:\bigwedge
V\rightarrow\bigwedge V$ be the unique odd superderivation such that
$(df)(v)=f(v)$ for any $v\in V\subseteq \bigwedge V$ (see, for
instance, \cite[\S 8]{KMRT}). Note that $(df)^2=0$.

Also, for any $v\in V$, the left multiplication by $v$ gives an odd
linear map $l_v:\bigwedge V\rightarrow\bigwedge V: x\mapsto vx$.
Again $l_v^2=0$, and for any $v\in V$ and $f\in V^*$:
\begin{equation}\label{eq:lvdf}
(l_v+df)^2=l_vdf+dfl_v=l_{(df)(v)}=f(v)id=q(v+f)id.
\end{equation}
The linear map $V\oplus V^*\rightarrow \End_k(\bigwedge V)$,
$v+f\mapsto l_v+df$, induces then an isomorphism
\begin{equation}\label{eq:Lambda}
\Lambda: \Cl(V\oplus V^*,q)\rightarrow \End_k(\bigwedge V).
\end{equation}
Moreover, let $\bar\ :\bigwedge V\rightarrow \bigwedge V$ be the
involution such that $\bar v=-v$ for any $v\in V$. Fix a basis
$\{v_1,\ldots,v_l\}$ of $V$, and let $\{f_1,\ldots,f_l\}$ be its
dual basis ($f_i(v_j)=\delta_{ij}$ for any $i,j=1,\ldots,l$). Let
$\Phi:\bigwedge V\rightarrow k$ be the linear function such that:
\begin{equation}\label{eq:Phi}
\Phi(v_1\cdots v_l)=1,\quad \Phi(v_{i_1}\cdots v_{i_r})=0\ \text{for
any $r<l$ and $1\leq i_1<\cdots<i_r\leq l$,}
\end{equation}
that is, $\Phi$ is a determinant, and consider the bilinear form
\begin{equation}\label{eq:b}
\begin{split}
b:\bigwedge V\times\bigwedge V&\rightarrow k\\
(s,t)\quad &\mapsto \Phi(\bar st).
\end{split}
\end{equation}
Since
\[
\begin{split}
\Phi(\overline{v_1\cdots v_l})&=(-1)^l\Phi(v_l\cdots v_1)=
 (-1)^l(-1)^{\binom{l}{2}}\Phi(v_1\cdots v_l)\\
 &=
 (-1)^{\binom{l+1}{2}}\Phi(v_1\cdots v_l)=(-1)^{\binom{l+1}{2}},
\end{split}
\]
it follows that, for any $s,t\in\bigwedge V$,
\[
b(t,s)=\Phi(\bar ts)=\Phi(\overline{\bar st})=
(-1)^{\binom{l+1}{2}}\Phi(\bar s t)=(-1)^{\binom{l+1}{2}}b(s,t).
\]
Hence,
\begin{equation}\label{eq:bsymskew}
\begin{split}
&\text{$b$ is symmetric if and only if $l\equiv 0$ or $3$
\hspace{-10pt}$\pmod
4$,}\\
&\text{$b$ is skew-symmetric if and only if $l\equiv 1$ or $2$
\hspace{-10pt}$\pmod 4$.}
\end{split}
\end{equation}

Let $\tau_b$ be the adjoint involution of $\bigwedge V$ relative to
$b$. Then, for any $v\in V$ and $s,t\in\bigwedge V$,
\[
b\bigl(l_v(s),t\bigr)=\Phi(\overline{vs}t)=\Phi(\bar s\bar vt)=
 -\Phi(\bar svt)=-b\bigl(s,l_vt\bigr),
\]
so $\tau_b(l_v)=-l_v$. Also, if $f\in V^*$ and $v\in V$,
\[
(df)(\bar v)=-(df)(v)=-f(v)=-\overline{f(v)}=(-1)^{\lvert
v\rvert}\overline{(df)(v)},
\]
where $\bigwedge V=\oplus_{i=0}^l \bigwedge^iV$ is the natural
$\bZ$-grading of $\bigwedge V$ and $\lvert s\rvert =i$ for
$s\in\bigwedge^iV$. Also, assuming $(df)(\bar s)=(-1)^{\lvert
s\rvert}\overline{(df)(s)}$ and $(df)(\bar t)=(-1)^{\lvert
t\rvert}\overline{(df)(t)}$ for homogeneous $s,t\in\bigwedge V$,
\[
\begin{split}
(df)(\overline{st})&=
 (df)(\bar t\bar s)=(df)(\bar t)\bar s +(-1)^{\lvert t\rvert}\bar
 t(df)(\bar s)\\
 &=(-1)^{\lvert t\rvert}\overline{(df)(t)}\bar s +
   (-1)^{\lvert s\rvert +\lvert t\rvert}\bar t\overline{(df)(s)}\\
 &=(-1)^{\lvert s\rvert +\lvert t\rvert}\overline{(df)(s)t+
      (-1)^{\lvert s\rvert}s(df)(t)}\\
 &=(-1)^{\lvert st\rvert}\overline{(df)(st)}.
\end{split}
\]
Hence $(df)(\bar s)=(-1)^{\lvert s\rvert}\overline{(df)(s)}$ for any
homogeneous $s\in\bigwedge V$. Thus, for any $f\in V^*$ and
$s,t\in\bigwedge V$,
\[
\begin{split}
b\bigl((df)(s),t\bigr)&=\Phi\bigl(\overline{(df)(s)}t\bigr)=
  (-1)^{\lvert s\rvert}\Phi\bigl((df)(\bar s)t\bigr)\\
  &=-\Phi\bigl(\bar s(df)(t)\bigr)\qquad\text{since
  $\Phi\bigl((df)(\bigwedge V)\bigr)=0$}\\
  &=-b\bigl(s,(df)(t)\bigr)
\end{split}
\]
and, therefore, $\tau_b(df)=-df$. As a consequence, the isomorphism
$\Lambda$ in \eqref{eq:Lambda} is actually an isomorphism of
algebras with involution:
\begin{equation}\label{eq:Lambdataus}
\Lambda:\bigl(\Cl(V\oplus V^*,q),\tau'\bigr)\rightarrow
\bigl(\End_k(\bigwedge V),\tau_b\bigr).
\end{equation}

\bigskip

The orthogonal Lie algebra $\frso_{2l+1}=\frso(W,q)$ is spanned by
the linear maps:
\begin{equation}\label{eq:sigmas}
\sigma_{w_1,w_2}=q(w_1,.)w_2-q(w_2,.)w_1
\end{equation}
where $q(w_1,w_2)=q(w_1+w_2)-q(w_1)-q(w_2)$ is the associated
symmetric bilinear form.

But for any $w_1,w_2,w_3\in W$, inside $\Cl(W,q)$ one has
\[
\begin{split}
[[w_1,w_2]^\cdot,w_3]^\cdot&=
 (w_1\cdot w_2-w_2\cdot w_1)\cdot w_3-w_3\cdot(w_1\cdot w_2-w_2\cdot
 w_1)\\
  &=q(w_2,w_3)w_1-w_1\cdot w_3\cdot w_2 -q(w_1,w_3)w_2+w_2\cdot w_3\cdot
  w_1\\
  &\quad -q(w_1,w_3)w_2+w_1\cdot w_3\cdot w_2+q(w_2,w_3)w_1-w_2\cdot w_3\cdot
  w_1\\
  &=-2\sigma_{w_1,w_2}(w_3).
\end{split}
\]
Therefore, $\frso_{2l+1}$ embeds in $\Cl\subo(W,q)$ by means of
$\sigma_{w_1,w_2}\mapsto -\frac{1}{2}[w_1,w_2]^\cdot$, so
$\frso_{2l+1}$ can be identified with the subspace $[W,W]^\cdot$ in
$\Cl\subo(W,q)$.

Under this identification, the action of $\frso_{2l+1}=\frso(W,q)$
on its natural module $W$ corresponds to the adjoint action of
$[W,W]^\cdot$ on $W$ inside $\Cl(W,q)$. Note that for any $x,y\in
V\oplus V^*$,
\[
\begin{split}
\Psi\bigl([x,y]^\cdot\bigr)=[u\cdot x,u\cdot y]^\cdot&=
  u\cdot x\cdot u\cdot y-u\cdot y\cdot u\cdot x\\
  &=-u\cdot u\cdot(x\cdot y-y\cdot x)\\
  &=x\cdot y-y\cdot x=[x,y]^\cdot,
\end{split}
\]
so $\Psi$ acts ``identically'' on $\frso_{2l}=[V\oplus V^*,V\oplus
V^*]^\cdot\subseteq \Cl(V\oplus V^*,q)$.

The subspace $\frh=\espan{[v_i,f_i]^\cdot: i=1,\ldots,l}$ is a
Cartan subalgebra of $\frso_{2l+1}\simeq [W,W]^\cdot$. Besides,
\[
\begin{split}
[[v_i,f_i]^\cdot,u]^\cdot&=0,\\
[[v_i,f_i]^\cdot,v_j]^\cdot&=2\delta_{ij}v_j,\\
[[v_i,f_i]^\cdot,f_j]^\cdot&=-2\delta_{ij}f_j.
\end{split}
\]
Hence, if $\epsilon_i:\frh\rightarrow k$ denotes the linear map with
$\epsilon_i\bigl([v_j,f_j]^\cdot\bigr)=2\delta_{ij}$, the weights of
the natural module $W$ relative to $\frh$ are $0$ and
$\pm\epsilon_i$, $i=1,\ldots,l$, all of them of multiplicity $1$;
while there appears a root space decomposition
\[
\frso_{2l+1}\simeq
[W,W]^\cdot=\frh\oplus\bigl(\oplus_{\alpha\in\Delta}\frg_\alpha\bigr),
\]
where
\[
\Delta=\{\pm(\epsilon_i+\epsilon_j): 1\leq i<j\leq l\}\cup
\{\pm\epsilon_i:1\leq i\leq
l\}\cup\{\pm(\epsilon_i-\epsilon_j):1\leq i<j\leq l\}.
\]
Here $\frg_{\epsilon_i+\epsilon_j}=k[v_i,v_j]^\cdot$,
$\frg_{-(\epsilon_i+\epsilon_j)}=k[f_i,f_j]^\cdot$,
$\frg_{\epsilon_i-\epsilon_j}=k[v_i,f_j]^\cdot$,
$\frg_{\epsilon_i}=k[u,v_i]^\cdot$, and
$\frg_{-\epsilon_i}=k[u,f_i]^\cdot$, for any $i\ne j$. This root
space decomposition induces a triangular decomposition
\[
\frso_{2l+1}=\frg_{-}\oplus\frh\oplus\frg_{+},
\]
where $\frg_{\pm}=\oplus_{\alpha\in\Delta^{\pm}}\frg_\alpha$, with
$\Delta^+=\{\epsilon_i+\epsilon_j: 1\leq i<j\leq l\}\cup
\{\epsilon_i:1\leq i\leq l\}\cup\{\epsilon_i-\epsilon_j:1\leq
i<j\leq l\}$, and $\Delta^-=-\Delta^+$.

\bigskip

The spin representation of $\frso_{2l+1}$ is given by the
composition
\[
\frso_{2l+1}\hookrightarrow
\Cl\subo(W,q)\overset{\Psi^{-1}}{\longrightarrow}\Cl(V\oplus
V^*,q)\overset{\Lambda}{\longrightarrow} \End_k(\bigwedge V).
\]
Denote this composition by $\rho$:
\begin{equation}\label{eq:rho}
\rho=\Lambda\circ\Psi^{-1}\vert_{\frso_{2l+1}},
\end{equation}
and denote by $S=\bigwedge V$ the spin module. Note that for any
$1\leq i\leq l$ and any $1\leq i_1<\cdots <i_r\leq l$
\[
\rho\bigl([v_i,f_i]^\cdot\bigr)(v_{i_1}\cdots v_{ir})=
\begin{cases} v_{i_1}\cdots v_{i_r}&\text{if $i=i_j$ for some
$j$},\\
 -v_{i_1}\cdots v_{i_r}&\text{otherwise}.
 \end{cases}
\]
Thus, $v_{i_1}\cdots v_{i_r}$ is a weight vector relative to $\frh$,
with weight
\[
\frac{1}{2}\bigl(\sum_{i\in\{i_1,\ldots,i_r\}}\epsilon_i-
\sum_{j\not\in\{i_1,\ldots,i_r\}}\epsilon_j\bigr),
\]
and hence all the weights of the spin module have multiplicity $1$.

\begin{proposition}\label{pr:spin}
Up to scalars, there is a unique $\frso_{2l+1}$-invariant bilinear
map $S\times S\rightarrow \frso_{2l+1}$, $(s,t)\mapsto [s,t]$. This
map is given by the formula
\begin{equation}\label{eq:spinproduct}
\frac{1}{2}\tr\bigl(\sigma[s,t]\bigr)=b\bigl(\rho(\sigma)(s),t\bigr),
\end{equation}
for any $\sigma\in\frso_{2l+1}$ and $s,t\in S$, where $\tr$ denotes
the trace of the natural representation of $\frso_{2l+1}$.

Moreover, this bilinear map $[.,.]$ is symmetric if and only if $l$
is congruent to $1$ or $2$ modulo $4$. Otherwise, it is
skew-symmetric.
\end{proposition}
\begin{proof}
First note that  the trace form $\tr$ is $\frso_{2l+1}$-invariant,
and so is $b$ because $\Psi$ and $\Lambda$ in \eqref{eq:Psitaus} and
\eqref{eq:Lambdataus} are isomorphisms of algebras with involutions.
Since both $\tr$ and $b$ are nondegenerate, $[.,.]$ is well defined
and $\frso_{2l+1}$-invariant. Now, the space of
$\frso_{2l+1}$-invariant bilinear maps $S\times S\rightarrow
\frso_{2l+1}$ is isomorphic to $\Hom_{\frso_{2l+1}}(S\otimes
S,\frso_{2l+1})$ (all the tensor products are considered over the
ground field $k$), or to the space of those tensors in $S\otimes
S\otimes (\frso_{2l+1})^*\simeq S\otimes S\otimes \frso_{2l+1}\simeq
\frso_{2l+1}\otimes  S\otimes S^*$ ($\tr$ and $b$ are nondegenerate)
annihilated by $\frso_{2l+1}$, and hence to
$\Hom_{\frso_{2l+1}}(\frso_{2l+1}\otimes S,S)$. But
$\frso_{2l+1}\otimes S$ is generated, as a module for
$\frso_{2l+1}$, by the tensor product of any nonzero element (like
$[v_1,v_2]^\cdot$) in the root space
$(\frso_{2l+1})_{\epsilon_1+\epsilon_2}$ ($(\frso_{3})_{\epsilon_1}$
if $l=1$), and any nonzero element (like $1$) in the weight space
$S_{-\frac{1}{2}(\epsilon_1+\cdots+\epsilon_l)}$. (Note that
$\epsilon_1+\epsilon_2$ is the longest root in the lexicographic
order given by $\epsilon_1>\cdots>\epsilon_l>0$, while
$-\frac{1}{2}(\epsilon_1+\cdots+\epsilon_l)$ is the lowest weight in
$S$.) The image of this basic tensor under any homomorphism of
$\frso_{2l+1}$-modules lies in the weight space of weight
$\frac{1}{2}(\epsilon_1+\epsilon_2-\epsilon_3-\cdots-\epsilon_l)$,
which is one dimensional. Hence,
$\dim_k\Hom_{\frso_{2l+1}}(\frso_{2l+1}\otimes S,S)=1$, as required.

The last part of the Proposition follows from \eqref{eq:bsymskew}.
\end{proof}

\medskip

For future use, note that for any $w_1,w_2,w_3,w_4\in W$,
\[
\tr\bigl(\sigma_{w_1,w_2}\sigma_{w_3,w_4}\bigr)=
2\bigl(q(w_1,w_4)q(w_2,w_3)-q(w_1,w_3)q(w_2,w_4)\bigr),
\]
and hence, under the identification $\frso_{2l+1}\simeq [W,W]^\cdot$
($\sigma_{w_1,w_2}\mapsto -\frac{1}{2}[w_1,w_2]^\cdot$),
\begin{equation}\label{eq:tr}
\frac{1}{2}\tr\bigl([w_1,w_2]^\cdot[w_3,w_4]^\cdot\bigr)=
4\bigl(q(w_1,w_4)q(w_2,w_3)-q(w_1,w_3)q(w_2,w_4)\bigr).
\end{equation}

\bigskip

In order to deal with the Lie algebras $\frso_{2l}$ (type $D$),
$l\geq 2$, consider the involution of $\Cl(V\oplus V^*,q)$, which
will be denoted by $\tau$ too,  which is the identity on $V\oplus
V^*$. Also consider the involution $\hat\ :\bigwedge V\rightarrow
\bigwedge V$ such that $\hat v=v$ for any $v\in V\subseteq \bigwedge
V$, and the nondegenerate bilinear form:
\begin{equation}\label{eq:bhat}
\begin{split}
\hat b: \bigwedge V\times\bigwedge V&\rightarrow k\\
 (s,t)\quad&\mapsto \Phi(\hat s t),
\end{split}
\end{equation}
where $\Phi$ is as in \eqref{eq:Phi}. Here, with the same arguments
as for \eqref{eq:bsymskew},
\begin{equation}\label{eq:bhatsymskew}
\begin{split}
&\text{$\hat b$ is symmetric if and only if $l\equiv 0$ or $1$
\hspace{-10pt}$\pmod 4$,}\\
&\text{$\hat b$ is skew-symmetric if and only if $l\equiv 2$ or $3$
\hspace{-10pt}$\pmod 4$.}
\end{split}
\end{equation}
Moreover, if $l$ is even, then $\hat b\bigl(\bigwedge\subo
V,\bigwedge\subuno V\bigr)=0$, so the restrictions of $\hat b$ to
$S^+=\bigwedge\subo V$ and $S^-=\bigwedge\subuno V$ are
nondegenerate. However, if $l$ is odd, then both $S^+$ and $S^-$ are
isotropic subspaces relative to $\hat b$.

The nondegenerate bilinear form $\hat b$ induces the adjoint
involution $\tau_{\hat b}$ on $\bigwedge V$ and, as before, the
isomorphism $\Lambda$ in \eqref{eq:Lambda} becomes an isomorphism of
algebras with involution:
\begin{equation}\label{eq:Lambdataushat}
\Lambda: \bigl(\Cl(V\oplus
V^*,q),\tau\bigr)\rightarrow\bigl(\End_k(\bigwedge V),\tau_{\hat
b}\bigr).
\end{equation}
Under this isomorphism, the even Clifford algebra $\Cl\subo(V\oplus
V^*,q)$ maps onto $\End_k(\bigwedge\subo
V)\oplus\End_k(\bigwedge\subuno V)$.

Also, as before, $\frso_{2l}=\frso(V\oplus V^*,q)$ can be identified
with the subspace $[V\oplus V^*,V\oplus V^*]^\cdot$ of
$\Cl\subo(V\oplus V^*,q)$, $\frh=\espan{[v_i,f_i]^\cdot:
i=1,\ldots,l}$ is a Cartan subalgebra, the roots are
$\{\pm\epsilon_i\pm\epsilon_j: 1\leq i<j\leq l\}$, the set of
weights of the natural module $V\oplus V^*$ are $\{\pm\epsilon_i:
1\leq i\leq l\}$, all the weights appear with multiplicity one, and
the composition
\[
\frso_{2l}\hookrightarrow \Cl\subo(V\oplus
V^*,q)\overset{\Lambda}{\longrightarrow}\End_k(\textstyle{\bigwedge\subo
V})\oplus\End_k(\textstyle{\bigwedge\subuno V})
\]
gives two representations
\begin{equation}\label{eq:rhopm}
\rho^+:\frso_{2l}\rightarrow \End_k(\textstyle{\bigwedge\subo
V})\quad\text{and}\quad \rho^-:\frso_{2l}\rightarrow
\End_k(\textstyle{\bigwedge\subuno V}),
\end{equation}
called the half-spin representations. The weights in
$S^+=\bigwedge\subo V$ (respectively $S^-=\bigwedge\subuno V$) are
the weights $\frac{1}{2}(\pm\epsilon_1\pm\cdots\pm\epsilon_l)$, with
an even (respectively odd) number of $+$ signs.

\begin{proposition}\label{pr:halfspin}
\null\quad\null
\begin{romanenumerate}
\item If $l$ is odd, $l\geq 3$, there is no nonzero
$\frso_{2l}$-invariant bilinear map $S^+\times S^+\rightarrow
\frso_{2l}$.
\item If $l$ is even, there is a unique, up to scalars, such
bilinear map, which is given by the formula
\begin{equation}\label{eq:halfspinproduct}
\frac{1}{2}\tr\bigl(\sigma[s,t]\bigr)=\hat
b\bigl(\rho^+(\sigma)(s),t\bigr),
\end{equation}
for any $\sigma\in\frso_{2l}$ and $s,t\in S^+$. Moreover, this
bilinear map $[.,.]$ is symmetric if and only if $l$ is congruent to
$2$ or $3$ modulo $4$, and it is skew-symmetric otherwise.
\end{romanenumerate}
\end{proposition}
\begin{proof}
If $l$ is odd ($l\geq 3$), then $S^+\otimes S^+$ is generated, as a
module for $\frso_{2l}$, by $v_1\cdots v_{l-1}\otimes v_{l-1}v_l$
(the tensor product or a nonzero highest weight vector and a nonzero
lowest weight vector), and its image under any nonzero
$\frso_{2l}$-invariant linear map $S^+\otimes S^+\rightarrow
\frso_{2l}$ lies in the root space of root
$\frac{1}{2}(\epsilon_1+\cdots+\epsilon_{l-1}-\epsilon_l)
+\frac{1}{2}(-\epsilon_1-\cdots-\epsilon_{l-2}
+\epsilon_{l-1}+\epsilon_l)=\epsilon_{l-1}$. But $\epsilon_{l-1}$ is
not a root, so its image must be $0$.

For $l$ even $\hat b$ is nondegenerate on $S^+$ and, as in
Proposition \ref{pr:spin}, it is enough to compute
$\dim_k\Hom_{\frso_{2l}}(\frso_{2l}\otimes S^+,S^+)$, which is
proven to be $1$ with the same arguments given there.
\end{proof}

\begin{remark}\label{rk:Sminus}
In $\Cl\subuno(V\oplus V^*,q)$ there are invertible elements $a$
such that $a^{\cdot 2}\in k1$ and $a\cdot(V\oplus V^*)\cdot
a^{-1}\subseteq V\oplus V^*$. For instance, one can take the element
$a=[v_1,f_1]^\cdot\cdots [v_{l-1},f_{l-1}]^\cdot\cdot (v_l+f_l)$,
which satisfies $a^{\cdot 2}=(-1)^{l-1}$. (Note that
$[v_i,f_i]^{\cdot 2}=v_i\cdot f_i\cdot v_i\cdot f_i+f_i\cdot
v_i\cdot f_i\cdot v_i=v_i\cdot(1-v_i\cdot f_i)\cdot f_i+f_i\cdot
(1-f_i\cdot v_i)\cdot v_i=v_i\cdot f_i+f_i\cdot v_i=f_i(v_i)=1$ and
$(v_i+f_i)\cdot(v_i+f_i)=v_i\cdot f_i+f_i\cdot v_i=1$.) Consider the
linear isomorphism $\phi_a:S^-=\bigwedge\subuno V\rightarrow
S^+=\bigwedge\subo V$, $s\mapsto \rho(a)(s)$, and the order two
automorphism $\Ad_a:\Cl(V\oplus V^*,q)\rightarrow \Cl(V\oplus
V^*,q)$, $x\mapsto a\cdot x\cdot a^{-1}$. $\Ad_a$ preserves $V\oplus
V^*$, and hence also $\frso_{2l}\simeq [V\oplus V^*,V\oplus
V^*]^\cdot$. Then, for any $\frso_{2l}$-invariant bilinear map
$[.,.]:S^+\times S^+\rightarrow \frso_{2l}$, one gets a
$\frso_{2l}$-invariant bilinear map $[.,.]^-:S^-\times
S^-\rightarrow \frso_{2l}$ at once by
\[
[s,t]^-=\Ad_a\bigl([\phi_a(s),\phi_a(t)]\bigr).
\]
Therefore, it is enough to deal with the half spin representation
$S^+$.
\end{remark}

\bigskip

\section{Type $B$}

Let $\frg=\frg\subo\oplus\frg\subuno$ be either a simple
$\bZ_2$-graded Lie algebra or a simple Lie superalgebra with
$\frg\subo=\frso_{2l+1}$ and $\frg\subuno=S$ (its spin module).

Because of Proposition \ref{pr:spin}, the product of two odd
elements can be assumed to be given by the bilinear map $[s,t]$ in
\eqref{eq:spinproduct}. Therefore, the possibilities for such a
$\frg$ are given precisely by the values of $l$ such that the
product $[s,t]$ satisfies the Jacobi identity:
\[
J(s_1,s_2,s_3)=\rho([s_1,s_2])(s_3)+\rho([s_2,s_3])(s_1)+\rho([s_3,s_1])(s_2)=0,
\]
for any $s_1,s_2,s_3\in S$. As in Section 2, $\rho$ denotes the spin
representation of $\frso_{2l+1}$.

But $S\otimes S\otimes S$ is generated, a a module for
$\frso_{2l+1}$, by the elements $1\otimes v_1\cdots v_l\otimes
v_{i_1}\cdots v_{i_r}$, where $0\leq r\leq l$ and $1\leq
i_1<\cdots<i_r\leq l$. As in Section 2, $\{v_1,\ldots,v_l\}$ denotes
a fixed basis of $V$ and $\{f_1,\ldots,f_l\}$ the corresponding dual
basis in $V^*$. The trilinear map $S\times S\times S\rightarrow S$,
$(s_1,s_2,s_3)\mapsto J(s_1,s_2,s_3)$ is $\frso_{2l+1}$-invariant,
so it is enough to check for which values of $l$ the Jacobian
\[
J(1,v_1\cdots v_l,v_{i_1}\cdots v_{i_r})
\]
is $0$ for any $0\leq r\leq l$ and $1\leq i_1<\cdots<i_r\leq l$. By
symmetry, it is enough to check the Jacobians
\[
J(1,,v_1\cdots v_l,v_1\cdots v_r)
\]
for $0\leq r\leq l$.

\begin{theorem}\label{th:typeB}
Let $l\in\bN$ and let $\frg=\frg\subo\oplus\frg\subuno$ be the
$\bZ_2$-graded algebra with $\frg\subo=\frso_{2l+1}$,
$\frg\subuno=S$ (its spin module), and multiplication given by the
Lie bracket of elements in $\frso_{2l+1}$, and by
\[
\begin{split}
&[\sigma,s]=-[s,\sigma]=\rho(\sigma)(s),\quad\text{$\rho$ as in
\eqref{eq:rho}},\\
&[s,t]\quad\text{given by \eqref{eq:spinproduct}}.
\end{split}
\]
for any $\sigma\in\frg\subo$ and $s,t\in\frg\subuno$. Then:
\begin{enumerate}
\item[\textup (i)] $\frg$ is a Lie algebra if and only if
either:\\[-8pt]
\begin{itemize}
\item $l=3$ and the characteristic of $k$ is $3$, and then $\frg$ is
isomorphic to the $29$ dimensional simple Lie algebra discovered by
Brown \cite{Brown82}, or
\item $l=4$, and then $\frg$ is isomorphic to the simple Lie algebra
of type $F_4$.\\[-8pt]
\end{itemize}
\item[\textup (ii)] $\frg$ is a Lie superalgebra if and only if
either:\\[-8pt]
\begin{itemize}
\item $l=1$, and then $\frg$ is isomorphic to the orthosymplectic
Lie superalgebra $\frosp(1,2)$, or
\item $l=2$, and then $\frg$ is isomorphic to the orthosymplectic
Lie superalgebra $\frosp(1,4)$, or
\item $l=5$ and the characteristic of $k$ is $5$, or
\item $l=6$ and the characteristic of $k$ is $3$.
\end{itemize}
\end{enumerate}
\end{theorem}
\begin{proof}
With the same notations as in Section 2, note that $v_1\cdots v_r$
is a weight vector relative to $\frh$, of weight
$\frac{1}{2}(\epsilon_1+\cdots+\epsilon_r-\epsilon_{r+1}-\cdots-\epsilon_l)$
for any $0\leq r\leq l$. Hence $[1,v_1\cdots
v_r]\in(\frso_{2l+1})_{-(\epsilon_{r+1}+\cdots+\epsilon_l)}$. In the
same vein, $[v_1\cdots v_l,v_1\cdots
v_r]\in(\frso_{2l+1})_{\epsilon_1+\cdots\epsilon_r}$. In particular,
\begin{equation}\label{eq:v1vr}
[1,v_1\cdots v_r]=0\ \text{if $0\leq r\leq l-3$};\quad [v_1\cdots
v_l,v_1\cdots v_r]=0\ \text{if $3\leq r\leq l$.}
\end{equation}
Also, $[1,v_1\cdots v_l]\in\frh$, so that $[1,v_1\cdots
v_l]=\sum_{i=1}^l\alpha_i[v_i,f_i]^\cdot$ for some
$\alpha_1,\ldots,\alpha_l\in k$. By \eqref{eq:spinproduct}
\begin{equation}\label{eq:alphai}
\frac{1}{2}\sum_{i=1}^l\alpha_i\tr\bigl(\sigma[v_i,f_i]^\cdot\bigr)
 =b\bigl(\rho(\sigma)(1),v_1\cdots v_l\bigr).
\end{equation}
Let $\sigma=[v_j,f_j]^\cdot$, then by \eqref{eq:tr}
\[
\frac{1}{2}\sum_{i=1}^l\alpha_i\tr\bigl([v_j,f_j]^\cdot
[v_i,f_i]^\cdot\bigr)
 =4f_i(v_j)f_j(v_i)=4\delta_{ij},
\]
while
$\rho\bigl([v_j,f_j]^\cdot\bigr)(1)=[l_{v_j},df_j](1)=-(df_j)(v_j)=-1$.
Thus, \eqref{eq:alphai} gives $4\alpha_j=-1$ for any $j=1,\ldots,l$,
so
\begin{equation}\label{eq:v1vl}
[1,v_1\cdots v_l]=-\frac{1}{4}\sum_{i=1}^l[v_i,f_i]^\cdot.
\end{equation}
In the same vein, $[1,v_1\cdots
v_{l-1}]\in(\frso_{2l+1})_{-\epsilon_l}$, so $[1,v_1\cdots
v_{l-1}]=\alpha[u,f_l]^\cdot$ for some $\alpha\in k$, and by
\eqref{eq:spinproduct}
\begin{equation}\label{eq:alpha}
\frac{1}{2}\alpha\tr\bigl([u,v_l]^\cdot[u,f_l]^\cdot\bigr)=
 b\bigl(\rho([u,v_l]^\cdot(1),v_1\cdots v_{l-1}\bigr).
\end{equation}
The left hand side is
$4\alpha\bigl(-q(u,u)q(v_l,f_l)\bigr)=8\alpha$, while in the right
hand side
$\rho\bigl([u,v_l]^\cdot\bigr)=2\rho\bigl(\Psi(v_l)\bigr)=2\Lambda(v_l)$,
so this side becomes
\[
\begin{split}
2b(v_l,v_1\cdots v_{l-1})&=2\Phi(\bar v_lv_1\cdots
v_{l-1})=-2\Phi(v_lv_1\cdots v_{l-1})\\
&=2(-1)^l\Phi(v_1\cdots v_l)=2(-1)^l.
\end{split}
\]
Therefore $\alpha=\frac{(-1)^l}{4}$ and
\begin{equation}\label{eq:v1vl1}
[1,v_1\cdots v_{l-1}]=\frac{(-1)^l}{4}[u,f_l]^\cdot.
\end{equation}
Similar arguments, which are left to the reader, give
\begin{align}
&[1,v_1\cdots v_{l-2}]=\frac{1}{2}[f_{l-1},f_l]^\cdot,\ \text{if
$l\geq 2$,}\label{eq:v1vl2}\\
&[v_1\cdots v_l,v_1]=
 -\frac{(-1)^{\binom{l+1}{2}}}{4}[u,v_1]^\cdot,\label{eq:v1vlv1}\\
&[v_1\cdots v_l,v_1v_2]=
 -\frac{(-1)^{\binom{l+1}{2}}}{2}[v_1,v_2]^\cdot.\label{eq:v1vlv1v2}
\end{align}
Now, if $l\geq 7$ and $3\leq r\leq l-3$,
\[
\begin{split}
J(1,v_1\cdots v_l,v_1\cdots v_r)&=
 [[1,v_1\cdots v_l],v_1\cdots v_r]\quad\text{(by \eqref{eq:v1vr})}\\
 &=-\frac{1}{4}\sum_{i=1}^l[[v_i,f_i]^\cdot,v_1\cdots v_r]
 \quad\text{(by \eqref{eq:v1vl})}\\
 &=-\frac{1}{4}\sum_{i=1}^l\rho\bigl([v_i,f_i]^\cdot\bigr)(v_1\cdots
 v_r)\\
 &=-\frac{1}{4}(r-(l-r))v_1\cdots v_r\\
 &=\frac{1}{4}(l-2r)v_1\cdots v_r.
\end{split}
\]
With $r=\frac{l-2}{2}$ if $l$ is even, or $\frac{l-1}{2}$ if $l$ is
odd, $l-2r\,(=1\ \text{or}\ 2)\ne 0$, so the Jacobi identity is not
satisfied.

Assume now that $l=6$, so $[s,t]$ is symmetric in $s,t\in S$ by
Proposition \ref{pr:spin}. Then
\[
\begin{split}
J(1,v_1\cdots v_6,1)&=2[[1,v_1\cdots v_6],1]\quad\text{(because
$[1,1]=0$ \eqref{eq:v1vr})}\\
 &=-\frac{1}{2}\sum_{i=1}^6[[v_i,f_i]^\cdot,1]\quad \text{(by
 \eqref{eq:v1vl})}\\
 &=-\frac{1}{2}\sum_{i=1}^6(-1)=3,
\end{split}
\]
so the characteristic of $k$ must be $3$. Assuming this is so, it is
easily checked that $J(1,v_1\cdots v_6,v_1\cdots v_r)=0$ for any
$0\leq r\leq 6$.

For $l=5$, the product $[s,t]$ is also symmetric (Proposition
\ref{pr:spin}) and, as before,
\[
J(1,v_1\cdots v_5,1)=2[[1,v_1\cdots
v_5],1]=-\frac{1}{2}\sum_{i=1}^5(-1)=\frac{5}{2},
\]
so the characteristic of $k$ must be $5$, and then $J(1,v_1\cdots
v_5,v_1\cdots v_r)=0$ for any $0\leq r\leq 5$.

For $l=4$, the product $[s,t]$ is skew-symmetric (Proposition
\ref{pr:spin}), so $J(1,v_1\cdots v_4,1)=J(1,v_1\cdots v_4,v_1\cdots
v_4)=0$ by skew-symmetry. The other instances of $J(1,v_1\cdots
v_4,v_1\cdots v_r)$, $1\leq r\leq 3$, are also checked to be
trivial.

With $l=3$, the product $[s,t]$ is skew-symmetric too. Hence, by
\eqref{eq:v1vl}, \eqref{eq:v1vl1}, and \eqref{eq:v1vl2},
\[
\begin{split}
J(1,v_1v_2v_3,v_1)&=
 [[1,v_1v_2v_3],v_1]+[[v_1v_2v_3,v_1],1]+[[v_1,1],v_1v_2v_3]\\
 &=-\frac{1}{4}\sum_{i=1}^3[[v_i,f_i]^\cdot,v_1]-
  \frac{1}{4}[[u,v_1]^\cdot,1]-\frac{1}{2}[[f_2,f_3]^\cdot,v_1v_2v_3]\\
 &= -\frac{1}{4}(1-1-1)v_1-\frac{1}{4}(2v_1)-\frac{1}{2}(-1-1)v_1=
  \frac{3}{4}v_1,
\end{split}
\]
and hence the characteristic must be $3$. The other instance of the
Jacobi identity to be checked: $J(1,v_1v_2v_3,v_1v_2)=0$, also holds
easily.

For $l=1$ or $l=2$, the Jacobi identity is satisfied too.

\smallskip

The assertions about which Lie algebras or superalgebras appear
follows at once, since all the algebras and superalgebras mentioned
in the statement of the Theorem satisfy the hypotheses. (For
$\frosp(1,4)$, the even part is isomorphic to the symplectic Lie
algebra $\frsp_4$, and the odd part is its natural $4$ dimensional
module. However $\frsp_4$ is isomorphic to $\frso_5$, and viewed
like this, the $4$ dimensional module is the spin module. The same
happens for $\frosp(1,2)$.)
\end{proof}

\begin{remark}
Up to our knowledge, the modular Lie superalgebras that occur for
$l=5$ and $l=6$ have not appeared previously in the literature. Note
that the simplicity of $\frso_{2l+1}$ and the irreducibility of its
spin module imply that these superalgebras are simple.
\end{remark}

\bigskip

\section{Type $D$}

In this section the situation in which $\frg\subo=\frso_{2l}$
($l\geq 2)$, and $\frg\subuno=S^+$ (half-spin module) will be
considered. First note that it does not matter which half-spin
representation is used (Remark \ref{rk:Sminus}). By Proposition
\ref{pr:halfspin}, it is enough to deal with even values of $l$.

\begin{theorem}\label{th:typeD}
Let $l$ be an even positive integer, and let
$\frg=\frg\subo\oplus\frg\subuno$ be the $\bZ_2$-graded algebra with
$\frg\subo=\frso_{2l}$, $\frg\subuno=S^+$, and multiplication given
by the Lie bracket of elements in $\frso_{2l}$, and by
\[
\begin{split}
&[\sigma,s]=-[s,\sigma]=\rho^+(\sigma)(s),\quad\text{$\rho^+$ as in
\eqref{eq:rhopm}},\\
&[s,t]\quad\text{given by \eqref{eq:halfspinproduct}}.
\end{split}
\]
for any $\sigma\in\frg\subo$ and $s,t\in\frg\subuno$. Then:
\begin{enumerate}
\item[\textup (i)] $\frg$ is a Lie algebra if and only if
either:\\[-8pt]
\begin{itemize}
\item $l=8$ and then $\frg$ is
isomorphic to the simple Lie algebra of type $E_8$, or
\item $l=4$, and then $\frg$ is isomorphic to the simple Lie algebra
$\frso_9$ (of type $B_4$).\\[-8pt]
\end{itemize}
\item[\textup (ii)] $\frg$ is a Lie superalgebra if and only if
either:\\[-8pt]
\begin{itemize}
\item $l=6$, and the characteristic of $k$ is $3$, and
then $\frg$ is isomorphic to the Lie superalgebra in \cite[Theorem
3.2(v)]{Eld05}, or
\item $l=2$, and then $\frg$ is isomorphic to the direct sum
$\frosp(1,2)\oplus\frsl_2$, of the orthosymplectic Lie superalgebra
$\frosp(1,2)$ and the three-dimensional simple Lie algebra.
\end{itemize}
\end{enumerate}
\end{theorem}
\begin{proof}
Note first that the restriction to $S^+=\bigwedge\subo V$ of the
bilinear form $\hat b$ in \eqref{eq:bhat} coincides with the
restriction of the bilinear form $b$ in \eqref{eq:b}. Then, as in
the proof of Theorem \ref{th:typeB}, the equations \eqref{eq:v1vr},
\eqref{eq:v1vl}, \eqref{eq:v1vl2}, and \eqref{eq:v1vlv1v2} are all
valid here.

If $l\geq 10$ and $4\leq r\leq l-4$, $r$ even,
\[
\begin{split}
J(1,v_1\cdots v_l,v_1\cdots v_r)&=
 [[1,v_1\cdots v_l],v_1\cdots v_r]\\
 &=-\frac{1}{4}\sum_{i=1}^l[[v_i,f_i]^\cdot,v_1\cdots v_r]\\
 &=-\frac{1}{4}(r-(l-r))v_1\cdots v_r\\
 &=\frac{1}{4}(l-2r)v_1\cdots v_r,
\end{split}
\]
so, with $r=\frac{l-2}{2}$ if $l$ is congruent to $2$ modulo $4$, or
$r=\frac{l-4}{4}$ otherwise, $l-2r$ equals $2$ or $4$, and the
Jacobi identity is not satisfied.

For $l=8$, $[s,t]$ is skew-symmetric and it is enough to check that
the Jacobian $J(1,v_1\cdots v_8,v_1\cdots v_r)$ is $0$ for $r=2$,
$4$ or $6$, which is straightforward.

For $l=6$, $[s,t]$ is symmetric and
\[
\begin{split}
J(1,v_1\cdots v_6,1)&=2[[1,v_1\cdots v_6],1]\\
&=-\frac{1}{2}\sum_{i=1}^6[[v_i,f_i]^\cdot,1]\\
&=-\frac{1}{2}\sum_{i=1}^6(-1)=3,
\end{split}
\]
so the characteristic of $k$ must be $3$ and then all the other
instances of the Jacobi identity hold.

For $l=4$, $[s,t]$ is skew-symmetric, and thus it is enough to deal
with $J(1,v_1v_2v_3v_4,v_1v_2)$:
\[
\begin{split}
J(1,&v_1v_2v_3v_4,v_1v_2)\\
&=
 [[1,v_1v_2v_3v_4],v_1v_2]+[[v_1v_2v_3v_4,v_1v_2],1]+
  [[v_1v_2,1],v_1v_2v_3v_4]\\
 &=-\frac{1}{4}\sum_{i=1}^4[[v_i,f_i]^\cdot,v_1v_2]-
  \frac{1}{2}[[v_1,v_2]^\cdot,1]-
  \frac{1}{2}[[f_1,f_2]^\cdot,v_1v_2v_3v_4]\\
 &=-\frac{1}{4}(1+1-1-1)v_1v_2-v_1v_2+v_1v_2=0.
\end{split}
\]
It is well-known that $\frg=\frso_9$ is $\bZ_2$-graded with
$\frg\subo=\frso_8$ and $\frg\subuno$ the natural module for
$\frso_8$. But here, the triality automorphism permutes the natural
and the two half-spin modules, so one can substitute the natural
module by any of its half-spin modules. Therefore, the Lie algebra
that appears is isomorphic to $\frso_9$.

Finally, for $l=2$, $[s,t]$ is symmetric and the Jacobi identity is
easily checked to hold. Since $\frso_4$ is isomorphic to
$\frsl_2\oplus\frsl_2$ and the two half-spin representations are the
two natural (two dimensional) modules for each of the two copies of
$\frsl_2$. It follows then that $\frg=\frg_1\oplus \frg_2$, where
$(\frg_1)\subo\simeq\frsl_2$ and $(\frg_1)\subuno$ the natural
module for $\frsl_2$ (and hence $\frg_1\simeq\frosp(1,2)$), while
$\frg_2=(\frg_2)\subo=\frsl_2$. Alternatively, the subspaces
$\espan{[v_1,v_2]^\cdot,[f_1,f_2]^\cdot,[v_1,f_1]^\cdot+[v_2,f_2]^\cdot,1,v_1v_2}$
and
$\espan{[v_1,f_2]^\cdot,[v_2,f_1]^\cdot[v_1,f_1]^\cdot-[v_2,f_2]^\cdot}$
are ideals of $\frg$, the first one being isomorphic to
$\frosp(1,2)$, and the second one to $\frsl_2$.
\end{proof}

\bigskip

\section{The Kac Jordan superalgebra and the Tits construction}

The aim of this section is to show that the Lie superalgebra in
Theorem \ref{th:typeB} for $l=5$ (and characteristic $5$) is related
to a well-known  construction by Tits, applied to the Cayley algebra
and the ten dimensional Kac Jordan superalgebra over $k$.

This last superalgebra is easily described in terms of the smaller
Kaplansky superalgebra \cite{BE}. The tiny Kaplansky superalgebra is
the three dimensional Jordan superalgebra $K=K\subo\oplus K\subuno$,
with $K\subo=ke$ and $K\subuno =U$, a two dimensional vector space
endowed with a nonzero alternating bilinear form $(.\vert .)$, and
multiplication given by
\[
e^2=e,\quad ex=xe=\frac{1}{2}x,\quad xy=(x\vert y)e,
\]
for any $x,y\in U$. The bilinear form $(.\vert .)$ can be extended
to a supersymmetric bilinear form by means of $(e\vert
e)=\frac{1}{2}$ and $(K\subo\vert K\subuno)=0$.

For any homogeneous $u,v\in K$, \cite[(1.61)]{BE} shows that
\begin{equation}\label{eq:LuLv}
[L_u,L_v]\,\bigl(=L_uL_v-(-1)^{\bar u\bar v}L_vL_u\bigr)=
\frac{1}{2}\bigl(u(v\vert .)-(-1)^{\bar u\bar v}v(u\vert .)\bigr),
\end{equation}
where $L_x$ denotes the left multiplication by $x$, $\bar x$ being
the degree of the homogeneous element $x$. Moreover, the Lie
superalgebra of derivations of $K$ is \cite{BE}:
\[
\der K=[L_K,L_K]=\frosp(K)\,\bigl(\simeq \frosp(1,2)\bigr).
\]

The Kac Jordan superalgebra is
\[
\cJ=k1\oplus (K\otimes K),
\]
with unit element $1$ and product determined \cite[(2.1)]{BE} by
\begin{equation}\label{eq:Kacproduct}
(a\otimes b)(c\otimes d)=(-1)^{\bar b\bar c}\Bigl(ac\otimes
bd-\frac{3}{4}(a\vert c)(b\vert d)1\Bigr),
\end{equation}
for homogeneous elements $a,b,c,d\in K$. Because of
\cite[Proposition 2.7 and Theorem 2.8]{BE}, the superspace spanned
by the associators $(x,y,z)=(xy)z-x(yz)=(-1)^{\bar y\bar
z}[L_x,L_z](y)$ is $(\cJ,\cJ,\cJ)=K\otimes K$, and the Lie
superalgebra of derivations of $\cJ$ is $\der
\cJ=[L_{\cJ},L_{\cJ}]$, which acts trivially on $1$ and leaves
invariant $(\cJ,\cJ,\cJ)=K\otimes K$. Considered then as subspaces
of $\End_k(K\otimes K)$
\begin{equation}\label{eq:derJ}
\der J=(\der K\otimes id)\,\oplus \,(id\otimes\der K)\,\bigl(\simeq
\frosp(1,2)\oplus\frosp(1,2)\bigr).
\end{equation}
More precisely \cite[(2.4)]{BE}, as endomorphisms of $K\otimes K$,
for any homogeneous $a,b,c,d\in K$,
\begin{equation}\label{eq:LabLcd}
[L_{a\otimes b},L_{c\otimes d}]= (-1)^{\bar b\bar c}\Bigl(
[L_a,L_c]\otimes (b\vert d)id\, +\, (a\vert c)id\otimes
[L_b,L_d]\Bigr).
\end{equation}
(It must be remarked here that, with the usual conventions for
superalgebras, $id\otimes\varphi$ acts on $a\otimes b$ as
$(-1)^{\bar\varphi\bar a}a\otimes\varphi(b)$ for homogeneous
$\varphi$ and $a,b$.) In particular,
\[
(\der\cJ)\subo=\bigl((\der K)\subo\otimes id\bigr)\, \oplus \,
\bigl(id\otimes (\der K)\subo)\bigr),
\]
and $(\der K)\subo$ is isomorphic to $\frsp(U)=\frsp_2=\frsl_2$
(acting trivially on the idempotent $e$). The restriction of
$(\der\cJ)\subo$ to the subspace $K\subuno\otimes K\subuno=U\otimes
U$ of $K\otimes K$ then gives an isomorphism
\[
(\der\cJ)\subo\cong \frso(U\otimes U)\,\bigl(=(\frsp(U)\otimes
id)\,\oplus \, (id\otimes\frsp(U))\bigr),
\]
where $U\otimes U$ is endowed with the symmetric bilinear form given
by
\[
(u_1\otimes u_2\vert v_1\otimes v_2)=(u_1\vert v_1)(u_2\vert v_2),
\]
for any $u_1,u_2,v_1,v_2\in U$.

\medskip

Assume now that the characteristic of $k$ is $\ne 2,3$. Let $(C,n)$
be a unital composition algebra over $k$ with norm $n$. That is, $n$
is a regular quadratic form satisfying $n(ab)=n(a)n(b)$ for any
$a,b\in C$. (See \cite[Chapter III]{Schafer} for the basic facts
about these algebras.)

Since the field $k$ is assumed to be algebraically closed, $C$ is
isomorphic to either $k$, $k\times k$, $\Mat_2(k)$ or the Cayley
algebra $\cC$ over $k$.

The map $D_{a,b}:C\rightarrow C$ given by
\begin{equation}\label{eq:Dab}
D_{a,b}(c)=[[a,b],c]-3(a,b,c)
\end{equation}
is the inner derivation determined by $a,b\in C$, and the Lie
algebra $\der C$ is spanned by these derivations. The subspace
$C^0=\{ a\in C: n(1,a)=0\}$ orthogonal to $1$ is invariant under
$\der C$.

For later use, let us state some properties of the inner derivations
of Cayley algebras. For any $a$, let $\ad_a=L_a-R_a$ ($L_a$ and
$R_a$ denote, respectively, the left and right multiplication by the
element $a$), and consider $\ad_C=\{\ad_a: a\in C\}=\{\ad_a:a\in
C^0\}$.

\begin{lemma}\label{le:C}
Let $\cC$ be the Cayley algebra over $k$ ($\charac k\ne 2,3$). Then,
\begin{romanenumerate}
\item $\cC^0$ is invariant under $\der C$ and $\ad_{\cC}$, both of
which annihilate $k1$. Moreover, as subspaces of $\End_k(\cC^0)$,
$\frso(\cC^0,n)=\der \cC\oplus\ad_{\cC}$.
\item $[\ad_a,\ad_b]=2D_{a,b}-\ad_{[a,b]}$ for any $a,b\in\cC$.
\item $D_{a,b}+\frac{1}{2}\ad_{[a,b]}=3\bigl(n(a,.)b-n(b,.)a\bigr)$
for any $a,b\in\cC^0$.
\end{romanenumerate}
\end{lemma}
\begin{proof}
First, in \cite[Chapter III, \S 8]{Schafer} it is proved that $\der
C$ leaves invariant $\cC^0$ and, as a subspace of $\End_k(\cC^0)$,
it is contained in the orthogonal Lie algebra $\frso(\cC^0,n)$. The
same happens for $\ad_{\cC}=\ad_{\cC^0}$, and
$\der\cC\cap\ad_{\cC}=0$. Hence, by dimension count,
$\frso(\cC^0,n)=\der\cC\oplus\ad_{\cC}$, which gives (i).

Now, $\cC$ is an alternative algebra. That is, the associator
$(a,b,c)=(ab)c-a(bc)$ is alternating on its arguments. Hence, for
any $a,b,c\in \cC$:
\[
\bigl(L_{ab}-L_aL_b\bigr)(c)=(a,b,c)=-(a,c,b)=[L_a,R_b](c).
\]
Interchange $a$ and $b$ and subtract to get
\[
L_{[a,b]}-[L_a,L_b]=2[L_a,R_b]
\]
and, similarly
\[
\begin{split}
&R_{[a,b]}+[R_a,R_b]=-2[L_a,R_b],\\
&\ad_{[a,b]}=[L_a,L_b]+[R_a,R_b]+4[L_a,R_b].
\end{split}
\]
Hence
\begin{equation}\label{eq:DabadabLaRb}
\begin{split}
D_{a,b}&=\ad_{[a,b]}-3(a,b,.)=\ad_{[a,b]}+3(a,.,b)\\
 &=\ad_{[a,b]}-3[L_a,R_b]\\
 &=[L_a,L_b]+[R_a,R_b]+[L_a,R_b],
\end{split}
\end{equation}
and
\[
\begin{split}
[\ad_a,\ad_b]&=[L_a-R_a,L_b-R_b]\\
 &=[L_a,L_b]+[R_a,R_b]-2[L_a,R_b]\\
 &=D_{a,b}-3[L_a,R_b]\\
 &=2D_{a,b}-\ad_{[a,b]},
\end{split}
\]
thus getting (ii).

Now, for any $a\in\cC$ (\cite[Chapter III, \S 4]{Schafer}):
\[
a^2-n(1,a)a+n(a)1=0,
\]
so for any $a\in\cC^0$, $a^2=-n(a)1$ and hence
\begin{equation}\label{eq:abbanab}
ab+ba=-n(a,b)1,
\end{equation}
for any $a,b\in\cC^0$. Therefore, for any $a,b,c\in\cC^0$:
\[
\begin{split}
2\bigl(n(a,c)b&-n(b,c)a\bigr)\\
 &=-(ac+ca)b-b(ac+ca)+(bc+cb)a-a(bc+cb)\\
 &=-(a,c,b)+(b,c,a)-(ca)b+(cb)a-b(ac)+a(bc)\\
 &=\bigl(2[L_a,R_b]+[R_a,R_b]+[L_a,L_b]\bigr)(c)\\
 &=\bigl(D_{a,b}+[L_a,R_b]\bigr)(c)\\
 &=\bigl(\frac{2}{3}D_{a,b}+\frac{1}{3}\ad_{[a,b]}\bigr)(c),
\end{split}
\]
because of \eqref{eq:DabadabLaRb}, which gives (iii).
\end{proof}

\bigskip

Let $J=J\subo\oplus J\subuno$ be now a unital Jordan superalgebra
with a normalized trace $t:J\rightarrow k$. That is, $t$ is a linear
map such that $t(1)=1$, and $t(J\subuno)=0=t\bigl((J,J,J)\bigr)$
(see \cite[\S 1]{BETits}). Then $J=k1\oplus J^0$, where $J^0=\{x\in
J: t(x)=0\}$, which contains $J\subuno$. For $x,y\in J^0$,
$xy=t(xy)1+x*y$, where $x*y=xy-t(xy)1$ is a supercommutative
multiplication on $J^0$. Since $(J,J,J)=[L_J,L_J](J)$ is contained
in $J^0$, the subspace $J^0$ is invariant under $\inder J=[L_J,L_J]$
(the Lie superalgebra of inner derivations).

\smallskip

Given a unital composition algebra $C$ and a unital Jordan
superalgebra with a normalized trace $J$, consider the superspace
\[
\cT(C,J)=\der C\oplus(C^0\otimes J^0)\oplus\inder J,
\]
with the superanticommutative product $[.,.]$ specified by (see
\cite{BETits}):

\begin{itemize}
\renewcommand{\itemsep}{4pt}
\item $\der C$ is a Lie subalgebra and $\inder J$ a Lie
subsuperalgebra of $\cT(C,J)$,

\item $[\der C,\inder J]=0$,

\item $[D,a\otimes x]=D(a)\otimes x$, $[d,a\otimes x]=a\otimes
d(x)$,

\item $[a\otimes x,b\otimes y]=t(xy)D_{a,b}+[a,b]\otimes
x*y-2n(a,b)[L_x,L_y]$,
\end{itemize}

\noindent for all $D\in\der C$, $d\in\inder J$, $a,b\in C^0$ and
$x,y\in J^0$.

If the Grassmann envelope $G(J)$ satisfies the Cayley-Hamilton
equation $ch_3(x)=0$ of $3\times 3$-matrices, where
\[
ch_3(x)=x^3-3t(x)x^2+\Bigl(\frac{9}{2}t(x)^2-\frac{3}{2}t(x^2)\Bigr)x-
 \Bigl(t(x^3)-\frac{9}{2}t(x^2)t(x)+\frac{9}{2}t(x)^3\Bigr)1,
\]
then $\cT(C,J)$ is known to be a Lie superalgebra (see
\cite[Sections 3 and 4]{BETits}).

This construction, for algebras, was considered by Tits
\cite{Tits66}, who used it to give a unified construction of the
exceptional simple Lie algebras. In the above terms, it was
considered in \cite{BZ} and \cite{BETits}.

\bigskip

The Kac Jordan superalgebra $\cJ$ is endowed with a unique
normalized trace, given necessarily by $t(1)=1$ and $t(K\otimes
K)=0$. Note that if $f=f^2$ is an idempotent  linearly independent
to $1$ in a unital Jordan superalgebra with a normalized trace $t$,
and if the Grassmann envelope satisfies the Cayley-Hamilton equation
$ch_3(x)=0$, in particular $ch_3(f)=0$ so, by linear independence,
$t(f)-\frac{9}{2}t(f)^2+\frac{9}{2}t(f)^3=0$, and $t(f)$ is either
$0$, $\frac{1}{3}$ or $\frac{2}{3}$. But, if $t(f)$ were $0$,
$0=ch_3(f)$ would be equal to $f^3\,(=f)$, a contradiction. Hence
$t(f)=\frac{1}{3}$ or $\frac{2}{3}$. In the Kac superalgebra $\cJ$,
the element $f=-\frac{1}{2}+2e\otimes e$ is an idempotent with
$t(f)=-\frac{1}{2}$. Hence the Grassmann envelope of $\cJ$ cannot
satisfy the Cayley-Hamilton equation of degree $3$ unless,
$-\frac{1}{2}=\frac{1}{3}$ or $-\frac{1}{2}=\frac{2}{3}$, that is,
unless the characteristic of $k$ be $5$ or $7$. Actually, McCrimmon
\cite{McC} has shown that the Grassmann envelope $G(\cJ)$ satisfies
this Cayley-Hamilton equation if and only if the characteristic is
$5$. In retrospect, this explains the appearance of the nine
dimensional pseudocomposition superalgebras over fields of
characteristic $5$ (and only over these fields) in \cite[Example 9,
Theorem 14 and concluding notes]{EOpseudo}.

\medskip

Assume from now on that the characteristic of the ground field $k$
is $5$.

Then, if $C$ is a unital composition algebra, then $\cT(C,\cJ)$ is
always a Lie superalgebra. Obviously $\cT(k,\cJ)=\inder \cJ=\der
\cJ$, which is isomorphic to $\frosp(1,2)\oplus\frosp(1,2)$ (see
\eqref{eq:derJ}), and $\cT(k\times k,\cJ)$ is naturally isomorphic
to $L_{\cJ^0}\oplus\der \cJ$ which, in turn, is isomorphic to
$\frosp(K\oplus K)\cong \frosp(2,4)$ (see \cite[Theorem 2.13]{BE}).
Also, it is well-known that $\cT(\Mat_2(k),\cJ)$ is isomorphic to
the Tits-Kantor-Koecher Lie superalgebra of $\cJ$, which is
isomorphic to the exceptional Lie superalgebra of type $F(4)$. (This
was used by Kac \cite{KacJordan} in his classification of the
complex finite dimensional simple Jordan superalgebras.)

Our final result shows that the Lie superalgebra $\cT(\cC,\cJ)$ is,
up to isomorphism, the simple Lie superalgebra in Theorem
\ref{th:typeB} for $l=5$.

\begin{theorem}
Let $\cC$ be the Cayley algebra and let $\cJ$ be the Kac Jordan
superalgebra over an algebraically closed field $k$ of
characteristic $5$. Then:
\begin{romanenumerate}
\item $\cT(\cC,\cJ)\subo$ is isomorphic to the orthogonal Lie
algebra $\frso_{11}$.
\item $\cT(\cC,\cJ)\subuno$ is isomorphic to the spin module for
$\cT(\cC,\cJ)\subo$.
\end{romanenumerate}
\end{theorem}
\begin{proof}
For (i) consider the vector space
\[
M=\cC^0\oplus(U\otimes U)
\]
(recall that $U=K\subuno$), endowed with the symmetric bilinear form
$Q$ such that
\[
\begin{split}
&Q(\cC^0,U\otimes U)=0,\\
&Q(x)=-n(x),\\
&Q(u_1\otimes u_2,v_1\otimes v_2)=-(u_1\vert
v_1)(u_2\vert v_2),
\end{split}
\]
for $x\in\cC^0$ and $u_1,u_2,v_1,v_2\in U$. It will be shown that
$\cT(\cC,\cJ)\subo$ is isomorphic to the orthogonal Lie algebra
$\frso(M,Q)$.

This last orthogonal Lie algebra is spanned by the maps
\[
\sigma_{x,y}^Q=Q(x,.)y-Q(y,.)x
\]
for $x,y\in M$, and for any $\sigma\in\frso(M,Q)$,
\[
[\sigma,\sigma_{x,y}^Q]=\sigma_{\sigma(x),y}^Q+\sigma_{x,\sigma(y)}^Q.
\]
Moreover, since $\cC^0$ and $U\otimes U$ are orthogonal relative to
$Q$,
\begin{equation}\label{eq:soMQ}
\frso(M,Q)=\Bigl(\frso\bigl(\cC^0,Q\vert_{\cC^0}\bigr)\oplus
\frso\bigl(U\otimes U, Q\vert_{U\otimes U}\bigr)\Bigr)\oplus
\sigma_{\cC^0,U\otimes U}^Q
\end{equation}
(which gives a $\bZ_2$-grading of $\frso(M,Q)$), where
$\frso\bigl(\cC^0,Q\vert_{\cC^0}\bigr)$ and $\frso\bigl(U\otimes U,
Q\vert_{U\otimes U}\bigr)$ are embedded in $\frso(M,Q)$ in a natural
way, and $\frso(M,Q)$ is generated, as a Lie algebra, by
$\sigma_{\cC^0,U\otimes U}^Q$. Besides, for any $a,b\in\cC^0$, and
$u_1,u_2,v_1,v_2\in U$,
\[
[\sigma_{a,u_1\otimes u_2}^Q,\sigma_{b,v_1\otimes v_2}^Q]
=Q(u_1\otimes u_2,v_1\otimes v_2)\sigma_{a,b}^Q +
Q(a,b)\sigma_{u_1\otimes u_2,v_1\otimes v_2}^Q.
\]
Also, by Lemma \ref{le:C}(iii), for any $a,b\in\cC^0$,
\begin{equation}\label{eq:sigmaabQ}
\sigma_{a,b}^Q=-2D_{a,b}-\ad_{[a,b]}.
\end{equation}

Now, the multiplication in $\cT(\cC,\cJ)$ gives, for any
$a,b\in\cC^0$, $u,u_1,u_2,v,v_1,$ $v_2\in U$, $D\in\der \cC$ and
$d\in(\der\cJ)\subo$:
\begin{subequations}\label{eq:TCJ}
\begin{align}
&[D,a\otimes(u\otimes v)]=D(a)\otimes (u\otimes v),\label{eq:TCJ1}\\
&[a\otimes(e\otimes e),b\otimes (u\otimes v)]=
 \frac{1}{4}[a,b]\otimes(u\otimes v)=-\ad_a(b)\otimes(u\otimes
 v),\label{eq:TCJ2}\\
&[d,a\otimes(u\otimes v)]=a\otimes d(u\otimes v),\label{eq:TCJ3}\\
&[D,d]=0=[d,a\otimes(e\otimes e)],\label{eq:TCJ4}\\
&[L_{u_1\otimes u_2},L_{v_1\otimes v_2}]\vert_{U\otimes U}=
  \frac{1}{2}\sigma_{u_1\otimes u_2,v_1\otimes v_2}^Q\vert_{U\otimes
  U},\label{eq:TCJ5}
\end{align}
\end{subequations}
since
\[
\begin{split}
(u_1\otimes u_2)\bigl((v_1\otimes v_2)(w_1\otimes w_2)\bigr)&=
Q(v_1\otimes v_2,w_1\otimes w_2)(u_1\otimes u_2)(e\otimes
e-\frac{3}{4}1)\\
 &=-\frac{1}{2}Q(v_1\otimes v_2,w_1\otimes w_2)(u_1\otimes u_2),
\end{split}
\]
for any $u_1,u_2,v_1,v_2,w_1,w_2\in U$.

Moreover, for any $a,b\in\cC^0$ and $u_1,u_2,v_1,v_2\in U$,
\[
\begin{split}
[a\otimes&(u_1\otimes u_2),b\otimes(v_1\otimes v_2)]\\
 &=t\bigl(u_1\otimes u_2)(v_1\otimes v_2)\bigr)D_{a,b}+
   [a,b]\otimes \bigl((u_1\otimes u_2)*(v_1\otimes v_2)\bigr)\\
   &\hspace{2in}
    -2n(a,b)[L_{u_1\otimes u_2},L_{v_1\otimes v_2}]\\
 &=-2Q(u_1\otimes u_2,v_1\otimes v_2)D_{a,b}
  +[a,b]\otimes\bigl(Q(u_1\otimes u_2,v_1\otimes v_2)(e\otimes
  e)\bigr)\\
  &\hspace{2in}
    -n(a,b)\bigl( 2[L_{u_1\otimes u_2},L_{v_1\otimes v_2}]\bigr)\\
 &=Q(u_1\otimes u_2,v_1\otimes
 v_2)\bigl(-2D_{a,b}+[a,b]\otimes(e\otimes e)\\
  &\hspace{2in} +Q(a,b)\bigl(
  2[L_{u_1\otimes u_2},L_{v_1\otimes v_2}]\bigr).
\end{split}
\]
Now, the equations in Lemma \ref{le:C} and equations
\eqref{eq:sigmaabQ} and \eqref{eq:TCJ} prove that the linear map
\[
\Phi\subo:\cT(\cC,\cJ)\subo=\der C\oplus\bigl(\cC^0\otimes(U\otimes
U)\bigr)\oplus(\der\cJ)\subo\rightarrow \frso(M,Q),
\]
such that
\begin{itemize}
\item $\Phi\subo(D)=D$ for any $D\in\der\cC\,\bigl(\subseteq
\frso(\cC^0,n)=\frso(\cC^0,-n)\subseteq \frso(M,Q)\bigr)$, for any
$D\in\der\cC$,
\item $\Phi\subo(d)=d\vert_{U\otimes U}\,\bigl(\in \frso(U\otimes
U,Q)\subseteq \frso(M,Q)\bigr)$, for any $d\in(\der\cJ)\subo$,
\item $\Phi\subo\bigl(a\otimes(e\otimes
e)\bigr)=-\ad_a\,\bigl(\in \frso(\cC^0,-n)\subseteq
\frso(M,Q)\bigr)$, for any $a\in\cC^0$,
\item $\Phi\subo\bigl(a\otimes(u\otimes v)\bigr)=\sigma_{a,u\otimes
v}^Q$, for any $a\in\cC^0$ and $u,v\in U$,
\end{itemize}
is an isomorphism of Lie algebras. This proves the first part of the
Theorem.

\medskip

Now, let us consider the linear map
\[
\begin{split}
\Psi: M&\longrightarrow \End_k\bigl(\cC\otimes (U\oplus U)\bigr)\\
a\in\cC^0&\mapsto\quad L_a\otimes\begin{pmatrix} -id&0\\
0&id\end{pmatrix},\\
u_1\otimes u_2&\mapsto id\otimes \begin{pmatrix} 0&(u_2\vert .)u_1\\
(u_1\vert .)u_2&0\end{pmatrix}.
\end{split}
\]
(The elements in $U\oplus U$ are written as $\binom{u_1}{u_2}$, and
then $\End_k(U\oplus U)$ is identified with
$\Mat_2\bigl(\End_k(U)\bigr)$.) For any $a\in\cC^0$ and
$u_1,u_2,v_1,v_2\in U$:
\[
\begin{split}
&\Psi(a)\Psi(u_1\otimes u_2)+\Psi(u_1\otimes u_2)\Psi(a)=0,\\[6pt]
&\Psi(a)^2=-n(a)id=Q(a)id\quad\text{(as $a(ab)=a^2b=-n(a)b$ since
$\cC$ is alternative),}\\[6pt]
&\Psi(u_1\otimes u_2)\Psi(v_1\otimes v_2)+
 \Psi(v_1\otimes v_2)\Psi(u_1\otimes u_2)\\
&\quad =id\otimes\begin{pmatrix} (u_2\vert v_2)(v_1\vert
.)u_1+(v_2\vert u_2)(u_1\vert .)v_1&0\\
0& (u_1\vert v_1)(v_2\vert .)u_2+(v_1\vert u_1)(u_2\vert .)v_2
\end{pmatrix}\\
&\quad = id\otimes\begin{pmatrix} -(u_1\vert v_1)(u_2\vert v_2)id&0\\
  0& -(u_1\vert v_1)(u_2\vert v_2)id\end{pmatrix}\\
&\quad =Q(u_1\otimes u_2,v_1\otimes v_2)id,
\end{split}
\]
since $(u_2\vert v_2)\bigl((v_1\vert w_1)u_1+(w_1\vert
u_1)v_1\bigr)=-(u_2\vert v_2)(u_1\vert v_1)w_1$, because $(u_1\vert
v_1)w_1+(v_1\vert w_1)u_1+(w_1\vert u_1)v_1=0$, as this is an
alternating trilinear map on a two dimensional vector space.

Therefore, $\Psi$ induces an algebra homomorphism
$\Cl(M,Q)\rightarrow \End_k\bigl(\cC\otimes(U\oplus U)\bigr)$, which
restricts to an isomorphism (by dimension count)
$\Psi:\Cl\subo(M,Q)\rightarrow \End_k\bigl(\cC\otimes(U\oplus
U)\bigr)$. Therefore, $\cC\otimes (U\oplus U)$ is the spin module
for $\frso(M,Q)$. Recall that $\frso(M,Q)$ embeds in $\Cl\subo(M.Q)$
by means of $\sigma_{x,y}^Q\mapsto -\frac{1}{2}[x,y]^\cdot$. Since
$\frso(M,Q)$ is generated by the elements $\sigma_{a,u_1\otimes
u_2}^Q$ ($a\in\cC^0$, $u_1,u_2\in U$), the spin representation is
determined by
\[
\begin{split}
\rho\bigl(\sigma_{a,u_1\otimes u_2}^Q\bigr)
 &=-\frac{1}{2}\Psi\bigl([a,u_1\otimes u_2]^\cdot\bigr)\\
 &=-\frac{1}{2}[\Psi(a),\Psi(u_1\otimes u_2)]\\
 &=-\Psi(a)\Psi(u_1\otimes u_2)\\
 &=L_a\otimes\begin{pmatrix} 0&(u_2\vert .)u_1\\ -(u_1\vert .)u_2&0
   \end{pmatrix}.
\end{split}
\]
Identify now $\cT(\cC,\cJ)\subo$ with $\frso(M,Q)$ through
$\Phi\subo$, and identify $\cT(\cC,\cJ)\subuno=\cC^0\otimes
\bigl((U\otimes e)\oplus (e\otimes U)\bigr)\,\oplus(\der\cJ)\subuno$
with $\cC\otimes(U\oplus U)$ by means of
\[
\begin{split}
\Phi\subuno:\cT(\cC,\cJ)\subuno&\longrightarrow \cC\otimes(U\oplus
U)\\
a\otimes(u_1\otimes e+e\otimes u_2)&\mapsto\quad
a\otimes\binom{u_1}{u_2},\\
[L_e,L_{u_1}]\otimes id + id\otimes [L_e,L_{u_2}]&\mapsto\
-\frac{1}{2}\left(1\otimes\binom{u_1}{u_2}\right),
\end{split}
\]
for $a\in\cC^0$ and $u_1,u_2\in U$.

In $\cT(\cC,\cJ)$, for any $a,b\in\cC^0$, $u_1,u_2,v_1,v_2\in U$,
using \eqref{eq:LabLcd} we get
\[
\begin{split}
[a\otimes(u_1&\otimes u_2),b\otimes(v_1\otimes e+e\otimes v_2)]\\
 &=[a,b]\otimes\frac{1}{2}\bigl(u_1\otimes (u_2\vert v_2)e-(u_1\vert
 v_1)e\otimes u_2\bigr)\\
 &\qquad -2n(a,b)[L_{u_1\otimes u_2},L_{v_1\otimes e+e\otimes
 v_2}]\\
 &=\frac{1}{2}[a,b]\otimes\bigl((u_2\vert v_2)u_1\otimes e -
    e\otimes (u_1\vert v_1)u_2\bigr)\\
    &\qquad -2n(a,b)\bigl(-\frac{1}{2}(u_1\vert
    v_1)id\otimes[L_{u_2},L_e]+\frac{1}{2}[L_{u_1},L_e]\otimes
    (u_2\vert v_2)id\bigr)\\
 &=\frac{1}{2}[a,b]\otimes\bigl((u_2\vert v_2)u_1\otimes e -
    e\otimes (u_1\vert v_1)u_2\bigr)\\
    &\qquad -\frac{1}{2}n(a,b)\Bigl(-2\bigl(
     [L_e,L_{(u_2\vert v_2)u_1}]\otimes id \, -\, id\otimes
    [L_e,L_{(u_1\vert v_1)v_2}]\bigr)\Bigr).
\end{split}
\]
That is,
\[
\bigl[ a\otimes(u_1\otimes u_2),\Phi\subuno^{-1}\Bigl(b\otimes
\binom{v_1}{v_2}\Bigr)\bigl]=
 \Phi\subuno^{-1}\Bigl(ab\otimes \begin{pmatrix}0&(u_2\vert .)u_1\\
 -(u_1\vert .)u_2&0\end{pmatrix}\binom{v_1}{v_2}\Bigr),
\]
or
\[
\bigl[\Phi\subo^{-1}\bigl(\sigma_{a,u_1\otimes u_2}^Q\bigr),
 \Phi\subuno^{-1}\Bigl(b\otimes
\binom{v_1}{v_2}\Bigr)\bigl]=
 \Phi\subuno^{-1}\Bigl(\rho\bigl(\sigma_{a,u_1\otimes u_2}^Q\bigr)
  \bigl(b\otimes\binom{v_1}{v_2}\Bigr),
\]
because $ab=-\frac{1}{2}n(a,b)+\frac{1}{2}[a,b]$ for any
$a,b\in\cC^0$ by \eqref{eq:abbanab}. But also,
\[
\begin{split}
\bigl[a\otimes(u_1\otimes u_2),&[L_e,L_{v_1}]\otimes id+
id\otimes[L_e,L_{v_2}]\bigr]\\
 &=a\otimes\Bigl(-[L_e,L_{v_1}](u_1)\otimes u_2+u_1\otimes
 [L_e,L_{v_2}](u_2)\Bigr)\\
 &=a\otimes\Bigl(\frac{1}{2}(u_1\vert v_1)e\otimes u_2
 -\frac{1}{2}u_1\otimes (u_2\vert v_2)e\Bigr),
\end{split}
\]
or
\[
\bigl[a\otimes(u_1\otimes u_2),\Phi\subuno^{-1}\bigl(1\otimes
\binom{v_1}{v_2}\bigr)\bigr]=\Phi\subuno^{-1}\Bigl(a\otimes
\begin{pmatrix} 0&(u_2\vert .)u_1\\ -(u_1\vert
.)u_2&0\end{pmatrix}\binom{v_1}{v_2}\Bigr),
\]
that is,
\[
\bigl[\Phi\subo^{-1}\bigl(\sigma_{a,u_1\otimes u_2}^Q\bigr),
\Phi\subuno^{-1}\bigl(1\otimes \binom{v_1}{v_2}\bigr)\bigr]=
 \Phi\subuno^{-1}\Bigl(\rho\bigl(\sigma_{a,u_1\otimes u_2}^Q\bigr)
 \bigl(1\otimes \binom{v_1}{v_2}\bigr)\Bigr),
\]
and this shows that, if $\cT(\cC,\cJ)\subo$ is identified with
$\frso(M,Q)$ through $\Phi\subo$ and $\cT(\cC,\cJ)\subuno$ with
$\cC\otimes (U\oplus U)$ through $\Phi\subuno$, then the action of
$\cT(\cC,\cJ)\subo$ on $\cT(\cC,\cJ)\subuno$ is given, precisely, by
the spin representation.
\end{proof}

\bigskip

The Lie superalgebra in Theorem \ref{th:typeB} for $l=6$ and
characteristic $3$, appears in the extended Freudenthal Magic Square
in this characteristic \cite{CE}, as the Lie superalgebra
$\frg\bigl(\overline{B(4,2)},\overline{B(4,2)}\bigr)$, associated to
two copies of the unique six dimensional symmetric composition
superalgebra. This is related to the six dimensional simple
alternative superalgebra $B(4,2)$ \cite{Shestakov}, and hence to the
exceptional Jordan superalgebra of $3\times 3$ hermitian matrices
over $B(4,2)$, which is exclusive of characteristic $3$.


\providecommand{\bysame}{\leavevmode\hbox
to3em{\hrulefill}\thinspace}
\providecommand{\MR}{\relax\ifhmode\unskip\space\fi MR }
\providecommand{\MRhref}[2]{%
  \href{http://www.ams.org/mathscinet-getitem?mr=#1}{#2}
} \providecommand{\href}[2]{#2}

\end{document}